\documentclass[a4paper,12pt]{article}
\usepackage{a4wide}
\usepackage{amsmath,amsfonts,amssymb,amsthm,url,textcomp,mathrsfs}
\usepackage{anyfontsize}
\usepackage{csquotes}
\usepackage{enumitem}
\usepackage{fancyhdr}
\usepackage[colorlinks=true, linkcolor=cyan]{hyperref}
\usepackage[capitalise]{cleveref}
\usepackage{leftidx}
\usepackage{mathtools}
\usepackage{mdframed}
\usepackage[numbers]{natbib}
\usepackage{xcolor}
\usepackage[ruled,vlined]{algorithm2e}

\hypersetup{citecolor=cyan}

\allowdisplaybreaks

\makeatletter
\let\@fnsymbol\@arabic
\makeatother

\numberwithin{cor}{subsection}
\newtheorem{Cor}{Corollary}\numberwithin{Cor}{section}

\newtheorem{dfn}[Cor]{Definition}

\newtheorem{exm}[Cor]{Example}

\newtheorem{lem}[Cor]{Lemma}

\newtheorem{prp}[Cor]{Proposition}

\newtheorem{rmk}[Cor]{Remark}

\newtheorem{thm}[Cor]{Theorem}

\newtheorem*{thm*}{Theorem}

\newcommand{\FE}[1]{\mathbb{F}_{2^{#1}}}
\newcommand{\FF}[2]{\mathbb{F}_{{#1}^{#2}}}
\newcommand{\FP}[1]{\mathbb{F}_{#1}}
\newcommand{\KK}{\mathbb{K}}

\newcommand{\NI}{\mathbb{N}_0}
\newcommand{\NN}{\mathbb{N}}
\newcommand{\QQ}{\mathbb{Q}}
\newcommand{\ZZ}{\mathbb{Z}}
\newcommand{\Tr}{\text{Tr}}
\newcommand{\Mod}[1]{\ (\mathrm{mod}\ #1)}

\begin{document}

\title{On stable polynomials of degrees \texorpdfstring{$2,3,4$}{}}

\date{}

\author{ Tong Lin\textsuperscript{1} \and Qiang Wang\thanks{School of Mathematics and Statistics, Carleton University, 1125 Colonel By Drive, Ottawa ON K1S 5B6, Canada. The authors were supported by the Natural Sciences and Engineering Research Council of Canada (RGPIN-2023-04673).
\newline \emph{E-mail addresses}:
\href{mailto:tonglin4@cmail.carleton.ca}{\color{cyan}{tonglin4@cmail.carleton.ca}}\ (T. Lin), 
\href{mailto:wang@math.carleton.ca}{\color{cyan}{wang@math.carleton.ca}}\ (Q. Wang).
}
}

\maketitle

\abstract{
Let $q$ be a prime power. We construct stable polynomials of the form $b^{m-1}(x+a)^m+c(x+a)+d$ over a finite field $\FP{q}$ for $m=2,3,4$ by Capelli's lemma.  When $m=3$ and $q$ is even,  we confirm the conjecture of Ahmadi and Monsef-Shokri \cite{AMS20a}  that the polynomial $f(x) = x^3 + x^2 + 1$ is stable over $\FP{2}$.  Moreover, 
 when $m=2$ and $q\equiv 1\Mod{4}$, we improve a lower bound of the number of quadratic stable polynomials by Gom\'ez-P\'erez and Nicol\'as~\cite{GPN10a}. }

\section{Introduction}\label{sec1}
Let $\KK$ be a field. A polynomial $t(x)\in \KK[x]$ is stable over $\KK$ if for each $n\in \NN$, the $n$-th iterate $t^{(n)}(x)=t(t(\dots t(t(x))))$ of $t$ is irreducible over $\KK$. 
Problems concerning stability of polynomials over fields date back to the $1980$s, when Odoni \cite[Proposition 4.1]{Odo85a}, motivated by a number-theoretic question involving prime divisors of a non-linear recurrence relation, discovered a stable polynomial $t_1(x)=x^2-x+1$ over $\QQ$.  On the other hand,   for any prime $p$ and integer $n\geq 3$,  Odoni proved in \cite[Corollary 1.6]{Odo85b} that the $n$-th iterate of $t_2(x)=x^p-x-1$ is reducible over the prime field $\FP{p}$. %These are one of the very first (counter-)examples of stable polynomials over a field. 
Since then,  stability of polynomials, especially those of low degrees, over various fields have been extensively studied.   For example,  Ahmadi et al. showed in \cite[Theorem 4, Corollary 11]{ALAS12a} that \emph{almost all} monic quadratic polynomials in $\ZZ[x]$ are stable over $\QQ$ and that no quadratic polynomial is stable over a finite field of characteristic $2$.   In $2010$, Gom\'ez-P\'erez and Nicol\'as first estimated the number of stable quadratic polynomials over a finite field of odd characteristic \cite[Theorem 1]{GPN10a}. Together with  with Ostafe and Sardonil,   they also estimated  the number of stable polynomials of an arbitrary degree $d\in \NN$ over a finite field of odd characteristic in 2014 \cite[Theorem 5.5]{GPNOS14a}. In $2012$, Jones and Boston gave necessary and sufficient conditions (in terms of the so-called adjusted critical orbits) for a quadratic polynomial to be stable over a finite field of odd characteristic \cite[Proposition 2.3]{JB12a} (and see \cite{JB20a} for errata). Based on their results, Ostafe and Shparlinski estimated the complexity of testing the stability of quadratic polynomials over a finite field of odd characteristic \cite[Corollary 2]{OS10a}. 

For polynomials of degree greater than $2$,  determining whether such a polynomial  is stable over a field remains challenging, and not many examples of stable polynomials of degree $\geq 3$ are known to date.  For example, according to \cite[Theorem 4.3]{GPNOS14a}, no polynomials of the form $ax^3+bx+c$ is stable over $\FF{3}{k}$. When it comes to cubic polynomials over $\FE{k}$, Ahmadi and Monsef-Shokri suggests that $f(x)=x^3+x^2+1$ is stable over $\FP{2}$ \cite[Conjecture 14]{AMS20a}.  Motivated by this conjecture,  we present several new classes of stable polynomials of degrees less than or equal to $4$ and thus confirm the conjecture of  Ahmadi and Monsef-Shokri.

The remaining part of this paper is organized as follows. In \cref{sec2}, we introduce a new class of polynomials, to which we apply the stability test proposed in \cite{{AMS20a}}. Based on the test, stable polynomials of the form
$$f(x)=b^{m-1}(x+a)^m+c(x+a)+d$$
where $m\in \{2,3,4\}$, are constructed in Sections \ref{sec3}, \ref{sec4} and \ref{sec5}, respectively. In addition, in \cref{sec3}, we improve the lower bound given in \cite[Theorem 1]{GPN10a} for the number of stable quadratic polynomials over $\FP{q}$ when $q\equiv 1\Mod{4}$. In \cref{sec4}, we verify the above-mentioned conjecture of  Ahmadi and Monsef-Shokri \cite[Conjecture 14]{AMS20a}, which is a special case of our construction. In \cref{sec5}, we also present an explicit result describing the factorization patterns for polynomials of the form $x^4+cx+d$ over $\FF{3}{m}$.  In \cref{sec6}, we make concluding remarks.

We need the following notations. Let $q$ be a prime power and $m\in \NN$. For each $d\in \NN$ such that $d\mid m$, the trace mapping from $\FF{q}{m}$ to $\FF{q}{d}$ is defined by

$$\begin{aligned}
\Tr_d^m(x) & =\sum_{j=0}^{\frac{m}{d}-1}x^{q^{jd}} \\ 
\end{aligned}$$

\noindent When $d=1$, we denote $\Tr_d^m$ by $\Tr_m$. 

\section{A sequence-based stability test}\label{sec2}
This section centers around a stability test based on the following Capelli's lemma.

\begin{lem}[\textbf{\cite[Lemma 13]{AMS20a}, \cite[Lemma 1]{Coh69a}}]\label{L21}
Let $q$ be any prime power, and let $f(x)\in \FP{q}[x]$ be an irreducible polynomial of degree $d\in \NN$. If $f^{(n)}(x)$ is irreducible over $\FP{q}$, then $f^{(n+1)}(x)$ is irreducible over $\FP{q}$ iff $f(x)-\alpha$ is irreducible over $\FF{q}{d^n}$ for some root $\alpha$ of $f^{(n)}(x)$ in $\FF{q}{d^n}$.
\end{lem}

\begin{rmk}\label{R22}
Capelli's Lemma provides a trade-off between polynomials of higher degrees (i.e., $\left(f^{(n)}(x)\right)_{n\geq 1}$) over a field of small size and polynomials of lower degrees (i.e., $\left(f(x)-\alpha_n\right)_{n\geq 1}$, where $\alpha_n\in \FF{q}{d^n}$ and $f^{(n)}(\alpha_n)=0$) over larger fields.
\end{rmk}

Next we define a new class of polynomials (which we call \textit{shift-resistant} polynomials) and study their stability by repeated use of Capelli's lemma.

\begin{dfn}\label{D23}
A polynomial $f(x)$ is called an $SR(q,d)$-polynomial (where $SR$ stands for shift-resistant) if $\deg(f)=d$ and for any $r\in \NN$ and any $a\in \FF{q}{r}$, $f(x)-a$ has a root in $\FF{q}{dr}$.
\end{dfn}

\begin{exm}\label{E24}
Every quadratic polynomial is an $SR(q,2)$-polynomial. Indeed, if $f(x)$ is quadratic, then $f(x)-a$ either has a root in $\FF{q}{r}$ or is irreducible over $\FF{q}{r}$. In the latter case, $f(x)-a$ has a root in its splitting field $\FF{q}{2r}$. Similarly, every cubic polynomial is an $SR(q,3)$-polynomial.
\end{exm}

\begin{exm}\label{E25}
Every quartic polynomial is an $SR(q,4)$-polynomial. Indeed, assume that $f(x)$ is quartic. If $f(x)-a$ has a root in $\FF{q}{r}$, then the claim is clearly true. Otherwise, there are two possibilities. If $f(x)-a$ is irreducible over $\FF{q}{r}$, then it has a root in its splitting field $\FF{q}{4r}$. If $f(x)-a$ is factored as a product of two irreducible quadratic factors over $\FF{q}{r}$, then it has a root in $\FF{q}{2r}\subset \FF{q}{4r}$.
\end{exm}

\begin{exm}\label{E26}
The quintic polynomial $f(x)=\left(x^2+x+1\right)\left(x^3+x+1\right)$ is not an $SR(2,5)$ polynomial because it has $2$ roots in $\FE{2}\setminus \FP{2}$ and $3$ roots in $\FE{3}\setminus \FP{2}$, and therefore none of its roots is in $\FE{5}$. 
\end{exm}

Let $f(x)$ be any $SR(q,d)$-polynomial. Let $\alpha_0=0$. For each $n\in \NN$, let $\alpha_n$ be a root of $f(x)-\alpha_{n-1}$ in $\FF{q}{d^n}$. Then it is easy to see that
$$f^{(n)}(\alpha_n)=f^{(n-1)}(f(\alpha_n))=f^{(n-1)}(\alpha_{n-1})=\dots=f(\alpha_1)=\alpha_0=0$$
In other words, $\alpha_n$ is a root of $f^{(n)}(x)$ in $\FF{q}{d^n}$. Then for each $n\in \NI=\NN\cup \{0\}$, let
$$f_n(x)=f(x)-\alpha_n$$
In the notation above, the following result is a restatement of Capelli's lemma.
\begin{prp}\label{P27}
Let $f(x)$ be any $SR(q,d)$-polynomial and let $n\in \NI$. Assume that $f^{(n)}(x)$ is irreducible over $\FP{q}$. Then $f^{(n+1)}(x)$ is irreducible over $\FP{q}$ iff $f_n(x)$ is irreducible over $\FF{q}{d^n}$.
\end{prp}

We note that $\alpha_0=0$ is a root of $f^{(0)}(x)=x$ in $\FP{q}$ and that $f_0(x)=f(x)$. So \cref{P27} holds trivially when $n=0$. Using the above result, a stability test for $SR(q,d)$-polynomials is obtained.

\begin{Cor}\label{C28}
Let $f(x)$ be any $SR(q,d)$-polynomial. Then $f(x)$ is stable over $\FP{q}$ iff for each $n\in \NI$, $f_n(x)$ is irreducible over $\FF{q}{d^n}$.
\end{Cor}

\noindent\textit{Proof}\quad If $f(x)$ is stable over $\FP{q}$, then $f^{(n)}(x)$ is irreducible over $\FP{q}$ for each $n\in \NI$. By \cref{P27}, $f_n(x)$ is irreducible over $\FF{q}{d^n}$ for each $n\in \NI$. Conversely, if $f_n(x)$ is irreducible over $\FF{q}{d^n}$ for each $n\in \NI$, then $f_0(x)=f(x)$ is irreducible over $\FP{q}$ and $f_1(x)$ is irreducible over $\FF{q}{d}$. By \cref{P27}, $f^{(2)}(x)$ is irreducible over $\FP{q}$. Then, as $f_2(x)$ is irreducible over $\FF{q}{d^2}$, $f^{(3)}(x)$ is irreducible over $\FP{q}$ due to \cref{P27}. In an iterative manner, we see that $f^{(n)}(x)$ is irreducible over $\FP{q}$ for each $n\in \NN$. Therefore, $f(x)$ is stable over $\FP{q}$.\mbox{}\hfill $\square$\\

In what follows, examples of stable $SR(q,d)$-polynomials for $d=2,3,4$ are given.

\section{Stable \texorpdfstring{$SR(q,2)$-\text{polynomials }($q\equiv 1\Mod{4}$)}{}}\label{sec3}
Let $q$ be an odd prime power. In this section, we study the stability of quadratic polynomials over $\FP{q}$ which are of the form
$$f(x)=b(x+a)^2+c(x+a)-a+\frac{c^2-2c}{4b}$$
where $b\neq 0$. 

\begin{thm}\label{T31}
If $a,b,c\in \FP{q}$, $b\neq 0$ and $4ab+2c$ is a non-square in $\FP{q}$, then
$$f(x)=b(x+a)^2+c(x+a)-a+\frac{c^2-2c}{4b}$$
is stable over $\FP{q}$ iff $q\equiv 1\Mod{4}$.
\end{thm}

\noindent\textit{Proof}\quad Following the idea introduced in \cref{sec2}, we construct a sequence $(\alpha_n)_{n\geq 0}$, where $\alpha_0=0$ and $f(\alpha_n)=\alpha_{n-1}\ (n\in \NN)$. For each $n\in \NI$, let
$$f_n(x)=f(x)-\alpha_n=b(x+a)^2+c(x+a)-(\alpha_n+a)+\frac{c^2-2c}{4b}$$
By \cref{C28}, $f(x)$ is stable over $\FP{q}$ iff $f_n(x)$ is irreducible over $\FF{q}{2^n}$ for each $n\in \NI$. The change of variables $y=b(x+a)$ implies that latter claim is true iff
$$g_n(y)=y^2+cy-\beta_n$$
is irreducible over $\FF{q}{2^n}$ for each $n\in \NI$, where
$$\beta_n=b(\alpha_n+a)-\frac{c^2-2c}{4}$$
Equivalently, $f(x)$ is stable over $\FP{q}$ iff
$$\Delta_n=c^2+4\beta_n$$
is a non-square in $\FF{q}{2^n}$ for each $n\in \NI$. It suffices to verify that the claim is true for all $n\in \NN$ because by assumption, $\Delta_0=4ab+2c$ is a non-square in $\FP{q}$. For each $n\in \NN$, since $f(\alpha_n)=\alpha_{n-1}$, we have that $g_{n-1}\left(\beta_n+\frac{c^2-2c}{4}\right)=0$, i.e., that
$$\left(\beta_n+\frac{c^2-2c}{4}\right)^2+c\left(\beta_n+\frac{c^2-2c}{4}\right)=\beta_{n-1}$$
meaning that
$$\left(\frac{\Delta_n-c^2}{4}+\frac{c^2-2c}{4}\right)^2+c\left(\frac{\Delta_n-c^2}{4}+\frac{c^2-2c}{4}\right)+\frac{c^2}{4}=\frac{\Delta_{n-1}-c^2}{4}+\frac{c^2}{4}$$
from which it follows that
\begin{equation}\label{Eq1}
\Delta_n^2=4\Delta_{n-1}
\end{equation}
\noindent \textbf{Case 1: }$q\equiv 3\Mod{4}$\\
\par In this case, $q+1\equiv 0\Mod{4}$. If $\gamma$ is primitive in $\FF{q}{2}$, then $\FP{q}=\left<\gamma^{q+1}\right>$. Hence, every element of $\FP{q}$ is a $4$-th power in $\FF{q}{2}$. So $\Delta_1^2\stackrel{(\ref{Eq1})}{=}4\Delta_0=s^4$ for some $s\in \FF{q}{2}$, meaning that $\Delta_1=(st)^2\in \FF{q}{2}$ for some $4$-th root of unity $t\in \FF{q}{2}$. Hence, $g_1(x)$ is reducible over $\FF{q}{2}$, and $f(x)$ is not stable over $\FP{q}$.\\
\\
\noindent \textbf{Case 2: }$q\equiv 1\Mod{4}$\\
\par In this case, we show that $\Delta_n$ a non-square in $\FF{q}{2^n}$ for each $n\in \NN$. Assume towards a contradiction that the claim is false. Then there exists a smallest $n_0\in \NN$ such that $\Delta_{n_0}$ is a square in $\FF{q}{2^{n_0}}$. If $n_0-1=0$, then $\Delta_{n_0-1}=\Delta_0$ is a non-square in $\FP{q}$. If $n_0-1\geq 1$, then by minimality of $n_0$, $\Delta_{n_0-1}$ is again a non-square in $\FF{q}{2^{n_0-1}}$. In either case, we have that 
$$\FF{q}{2^{n_0}}\cong \FF{q}{2^{n_0-1}}[x]\slash \left<x^2-\Delta_{n_0-1}\right>\cong \left\{u+vx_0:\ u,v\in \FF{q}{2^{n_0-1}}\right\}$$
where $x_0\in \FF{q}{2^{n_0}}$ and $x_0^2=\Delta_{n_0-1}$. Since $\Delta_{n_0}$ is a square in $\FF{q}{2^{n_0}}$, $\Delta_{n_0}=(u+vx_0)^2$ for some $u,v\in \FF{q}{2^{n_0-1}}$. Since $4x_0^2=4\Delta_{n_0-1}\stackrel{(\ref{Eq1})}{=}\Delta_{n_0}^2$, there is an $\epsilon\in \{\pm 1\}$ such that
$$2x_0=\epsilon\Delta_{n_0}=\epsilon(u+vx_0)^2=\epsilon\left(u^2+v^2\Delta_{n_0-1}\right)+2\epsilon uvx_0$$
which indicates that
$$\begin{cases}
\epsilon\left(u^2+v^2\Delta_{n_0-1}\right) & =0 \\
2\epsilon uv                               & =2
\end{cases}$$
The second equation implies that $v\neq 0$, which, together with the first equation, shows that $\Delta_{n_0-1}=-u^2v^{-2}$. Since $q\equiv 1\Mod{4}$, $-1$ is a square in $\FP{q}$. Thus, $\Delta_{n_0-1}$ is a square in $\FF{q}{2^{n_0-1}}$, which contradicts to the minimality of $n_0$. Hence, no such $n_0$ exists, proving that $f(x)$ is stable over $\FP{q}$.\mbox{}\hfill $\square$\\

\par In \cite[Theorem 1]{GPN10a}, Gom\'ez-P\'erez and Nicol\'as showed that for every odd prime power $q$, the number of (possibly non-monic) stable quadratic polynomials over $\FP{q}$, is at least $\frac{(q-1)^2}{4}$. When $q\equiv 1\Mod{4}$, we can improve this lower bound. Let $f(x)$ be as in \cref{T31}. We note that if $\delta=4ab+2c$, then
\begin{equation}\label{Eq2}
f(x)=bx^2+\left(\frac{\delta}{2}\right)x+\frac{\delta^2-4\delta}{16b}
\end{equation}
and that the discriminant of $f(x)$ is $\delta$. So $f(x)$ depends only on $b$ (for which there are $q-1$ choices as $b\neq 0$) and $\delta$ (which must be one of the $\frac{q+1}{2}$ non-squares in $\FP{q}$ if we want $f(x)$ to be stable over $\FP{q}$). These observations lead to the following result.

\begin{Cor}\label{C32}
If $q$ is a prime power and $q\equiv 1\Mod 4$, then
$$S_q\geq \frac{q^2-1}{2}$$
where $S_q$ is the number of (possibly non-monic) stable quadratic polynomials over $\FP{q}$.
\end{Cor}

The reason why we didn't rewrite $f(x)$ as \cref{Eq2} in the proof of \cref{T31} is that the change of variables $y=b(x+a)$ wouldn't be as obvious if we did. However, \cref{T31} can now be restated as follows.

\begin{thm}\label{T33}
Let $q$ be a prime power such that $q\equiv 1\Mod{4}$. If $b,\eta\in \FP{q}^\ast$ and
$$f(x)=bx^2+\left(\frac{\delta}{2}\right)x+\frac{\delta^2-4\delta}{16b}$$
then $f(x)$ is stable over $\FP{q}$ iff $\delta$ is a non-square in $\FP{q}$. Moreover, there are $\frac{q^2-1}{2}$ stable polynomials of the above-mentioned form.
\end{thm}

\section{Stable \texorpdfstring{$SR\left(2^m,3\right)$-\text{polynomials}}{}}\label{sec4}
Let $m\in \NN$. In this section, we study the stability of cubic polynomials over $\FE{m}$ which are of the form
$$f(x)=b^2(x+a)^2+x$$
where $a,b\neq 0$. In order to do so, the following results are needed.

\begin{prp}[\textbf{\cite[Theorem 7.20]{Wan03a}}]\label{P41}
If $m\in \NN$, $a,b,c\in \FE{m}$ and $a\neq 0$, then in $\FE{m}$, the equation $ax^2+bx+c=0$ has

$$\begin{cases}
\begin{aligned}
&\text{a unique solution}, & &\text{ if }\ b=0;\\
&\text{two distinct solutions}, & &\text{ if }\ b\neq 0\text{ and }\Tr_m\left(acb^{-2}\right)=0;\\
&\text{no solution},& &\text{ if }\ b\neq 0\text{ and }\Tr_m\left(acb^{-2}\right)=1.
\end{aligned}
\end{cases}$$
\end{prp}

\begin{lem}\label{L42}
If $m\in \NN, u,v\in \FE{m}^\ast$ and $u^3+u=v$, then $\Tr_m\left(u^{-1}\right)=\Tr_m\left(v^{-1}\right)$. 
\end{lem}

\noindent\textit{Proof}\quad We note that $x^2+(u+1)x+u=(x+u)(x+1)$ and $u\neq 1$ (else $v=0$). So
\begin{equation}\label{Eq3}
\Tr_m\left(u\left(u^2+1\right)^{-1}\right)
=\Tr_m\left(u\left(u+1\right)^{-2}\right)
=0
\end{equation}
where the last step follows from \cref{P41}. Since $v=u^3+u=u\left(u^2+1\right)$,
\begin{equation}\label{Eq4}
u^{-1}+u\left(u^2+1\right)^{-1}=\left(u\left(u^2+1\right)\right)^{-1}\left(\left(u^2+1\right)+u^2\right)=v^{-1}
\end{equation}
Applying $\Tr_m$ to \cref{Eq4} yields that
$$\begin{aligned}
\Tr_m\left(u^{-1}\right)\
&\stackrel{\mathmakebox[\widthof{=}]{(\ref{Eq3})}}{=}\ \Tr_m\left(u^{-1}\right)+\Tr_m\left(u\left(u^2+1\right)^{-1}\right) \\
&=\ \Tr_m\left(u^{-1}+u\left(u^2+1\right)^{-1}\right) \\
&\stackrel{\mathmakebox[\widthof{=}]{(\ref{Eq4})}}{=}\ \Tr_m\left(v^{-1}\right)
\end{aligned}$$
as required.\mbox{}\hfill $\square$\\

\par The theorem shown below, which concerns the factorization patterns of cubic trinomials over $\FE{m}$, plays a key role in proving our main result.

\begin{thm}[\textbf{\cite[Theorem 1]{KW75a}}]\label{T43}
Let $m\in \NN$. If $a,b\in \FE{m}$, $b\neq 0$ and $t_1,t_2$ are the two roots of $x^2+bx+a^3$, then over $\FE{m}$, the polynomial $x^3+ax+b$

\begin{enumerate}[label = (\arabic*)]
\item splits completely into linear factors iff $\Tr_m\left(a^3b^{-2}\right)=\Tr_m(1)$ and $t_1,t_2$ are both cubes in $\FE{m}$ (if $m$ is even) or $\FE{2m}$ (if $m$ is odd);
\item splits into a linear factor and a quadratic factor iff $\Tr_m\left(a^3b^{-2}\right)\neq \Tr_m(1)$;
\item is irreducible iff $\Tr_m\left(a^3b^{-2}\right)=\Tr_m(1)$ and $t_1,t_2$ are both non-cubes in $\FE{m}$ (if $m$ is even) or $\FE{2m}$ (if $m$ is odd).
\end{enumerate}
\end{thm}

Having done all the preparatory work, we now state the main result of this section.

\begin{thm}\label{T44}
If $m\in \NN$, $a,b\in \FE{m}^\ast$ and $\Tr_m\left((ab)^{-1}\right)=\Tr_m(1)$, then
$$f(x)=b^2(x+a)^3+x$$
is stable over $\FE{m}$ iff it is irreducible over $\FE{m}$.
\end{thm}

\noindent \textit{Proof}\quad We only prove that irreducibility implies stability since the other direction is trivial. Without loss of generality, we may consider only the case when $m$ is even as the other case can be proven in the same way. Let $q=2^m$. Following the idea introduced in \cref{sec2}, we construct a sequence $(\alpha_n)_{n\geq 0}$ such that $\alpha_0=0$ and that $f(\alpha_n)=\alpha_{n-1}\ (n\in \NN)$. For each $n\in \NI$, let
$$f_n(x)=f(x)+\alpha_n=b^2(x+a)^3+(x+a)+(\alpha_n+a)$$
By \cref{C28}, $f(x)$ is stable over $\FP{q}$ iff $f_n(x)$ is irreducible over $\FF{q}{3^n}$ for each $n\in \NI$. The change of variables $y=b(x+a)$ implies that latter claim is true iff
$$g_n(y)=y^3+y+\beta_n$$
is irreducible over $\FF{q}{3^n}$ for each $n\in \NI$, where
$$\beta_n=b(\alpha_n+a)$$
Since $f(\alpha_n)=\alpha_{n-1}$ for each $n\in \NN$, we have that $g_{n-1}(\beta_n)=0$, i.e., that
\begin{equation}\label{Eq5}
\beta_n^3+\beta_n=\beta_{n-1}
\end{equation}
Since $\beta_0=ab\neq 0$, it is easy to verify using \cref{Eq5} that $\beta_n\not\in 0,1$ for each $n\in \NN$. Moreover, \cref{Eq5} shows that any field containing $\beta_n$ also contains $\beta_j\ (0\leq j\leq n-1)$. These, together with repeated use of \cref{L42}, imply that
\begin{equation}\label{Eq6}
\Tr_{3^nm}\left(\beta_n^{-1}\right)=\Tr_{3^nm}\left(\beta_{n-1}^{-1}\right)=\dots=\Tr_{3^nm}\left(\beta_0^{-1}\right)
\end{equation}
Then we note that
$$\Tr_{3^nm}\left(\beta_0^{-1}\right)=\Tr_m\left(\Tr^{3^nm}_m\left((ab)^{-1}\right)\right)=\Tr_m\left((ab)^{-1}\right)=\Tr_m(1)=\Tr_{3^nm}(1)$$
where the last step is true as $3^nm\equiv m\Mod{2}$. We then deduce from \cref{Eq6} that
\begin{equation}\label{Eq7}
\Tr_{3^nm}\left(\beta_n^{-1}\right)=\Tr_{3^nm}(1)
\end{equation}

Since $m$ is even, \cref{Eq7} and \cref{P41} imply that $x^2+\beta_nx+1$ always has a root in $\FF{q}{3^n}=\FE{3^nm}$.

Assume towards a contradiction that $f(x)$ is not stable over $\FE{m}$. Then $g_n(x)$ is reducible over $\FF{q}{3^n}$ for some $n\in \NI$. Since $f_0(x)=f(x)$ is, by assumption, irreducible over $\FP{q}$, we know that $g_0(x)$ is irreducible over $\FP{q}$. Hence, there exists a smallest $n_0\in \NN$ such that $g_{n_0}(x)$ is reducible over $\FF{q}{3^{n_0}}$.

Applying \cref{T43} to $g_{n_0}(x)$, we deduce from \cref{Eq7} that the two roots $t_1,t_2$ of $x^2+\beta_{n_0}x+1$ are both cubes in $\FF{q}{3^{n_0}}$. If $t_1=w^3$ for some $w\in \FF{q}{3^{n_0}}^\ast$, then by Vieta's formulas, $t_2=t_1^{-1}=w^{-3}$ and
\begin{equation}\label{Eq8}
w^3+w^{-3}=t_1+t_2=\beta_{n_0}
\end{equation}
Let $u$ be a root of $x^2+\beta_{n_0-1}x+1$ in $\FF{q}{3^{n_0-1}}$. By Vieta's formula, $u^{-1}$ is the other root. Furthermore, we know that
\begin{equation}\label{Eq9}
w^9+w^{-9}=\left(w^3+w^{-3}\right)^3+\left(w^3+w^{-3}\right)\stackrel{(\ref{Eq8})}{=}\beta_{n_0}^3+\beta_{n_0}\stackrel{(\ref{Eq5})}{=}\beta_{n_0-1}
\end{equation}
Multiplying both sides of \cref{Eq9} by $w^9$, we see that $w^9$ is a root of $x^2+\beta_{n_0-1}x+1$. So $w^9=u$ or $u^{-1}$. In either case, $w^9\in \FF{q}{3^{n_0-1}}$. Letting $d=q^{3^{n_0-1}}$, we have that
\begin{equation}\label{Eq10}
\begin{aligned}
w & \in \FF{q}{3^{n_0}}^\ast=\FP{d^3}^\ast &\ \Rightarrow\ &\ w^{d^3-1}=1 &\ \Rightarrow\ &\ w^{\left(d^2+d+1\right)(d-1)}=1 \\
w^9 & \in \FF{q}{3^{n_0-1}}^\ast=\FP{d}^\ast &\ \Rightarrow\ &\ \left(w^9\right)^{(d-1)}=1 &\ \Rightarrow\ &\ \left(w^3\right)^{3(d-1)}=1
\end{aligned}
\end{equation}
Then we observe that $d^2+d+1\equiv 0\Mod{3}$ because when $m$ is even,
$$d\equiv \left(2^m\right)^{3^{n_0-1}}\equiv \left((-1)^m\right)^{3^{n_0-1}}\equiv 1\Mod{3}$$
We also note that $j^2+j+1\not\equiv 0\Mod{9}$ for each $j\in \NN$. So $d^2+d+1=3r$ for some $r\in \NN$ such that $3\nmid r$. By \cref{Eq10}, $\left(w^3\right)^{r(d-1)}=w^{\left(d^2+d+1\right)(d-1)}=1$, meaning that the order of $w^3$ in $\FP{d^3}^\ast$ divides $\gcd(r(d-1),3(d-1))=d-1$. In particular, $w^3\in \FP{d}^\ast=\FF{q}{3^{n_0-1}}^\ast$. Since either $u=\left(w^3\right)^3, u^{-1}=\left(w^{-3}\right)^3$ or the other way around, the two roots $u,u^{-1}$ of $x^2+\beta_{n_0-1}x+1$ are both cubes in $\FF{q}{3^{n_0-1}}$. By \cref{Eq7} and \cref{T43}, $g_{n_0-1}(x)$ is reducible over $\FF{q}{3^{n_0-1}}$. Since $n_0$ is the smallest positive integer such that $g_{n_0}(x)$ is reducible over $\FF{q}{3^{n_0}}$, we have that $n_0-1=0$. However, as has been shown earlier, $g_0(x)$ is irreducible over $\FP{q}$. Therefore, no such $n_0$ exists, giving the result.\mbox{}\hfill $\square$\\

\par When $a=b=1$, the exact same proof leads to the following result.

\begin{Cor}\label{C45}
The polynomial $f(x)=x^3+x+1$ is stable over $\FE{m}$ iff $3\nmid m$.
\end{Cor}

In \cref{C45}, when $m=1$, we obtain \cite[Conjecture 14]{AMS20a}.

\section{Stable \texorpdfstring{$SR(q,4)$-\text{polynomials for odd }$q$}{}}\label{sec5}
In this section, we give examples of stable quartic polynomials over $\FP{q}$ where $q$ is an odd prime power. To begin with, we recall the following result. 

\begin{thm}[\textbf{\cite[Theorem 8]{AMS20a}}]\label{T51}
Let $q$ be an odd prime power and let $f(x)$ be a polynomial of degree $d$ which is a product of $r$ pairwise distinct irreducible polynomials over $\FP{q}$. Then $r\equiv d\Mod{2}$ iff the discriminant of $f$ is a square in $\FP{q}$.
\end{thm}

The following is an immediate consequence of \cref{T51}.

\begin{prp}\label{P52}
If $q$ is a prime power and $q\equiv 1\Mod{4}$, then $f(x)=x^4-a$ is irreducible over $\FP{q}$ iff $a$ is a non-square in $\FP{q}$.
\end{prp}

\noindent\textit{Proof}\quad Since $4\mid q-1$, $-1$ is a square in $\FP{q}$. If $f(x)$ is irreducible over $\FP{q}$, then by \cref{T51}, the discriminant of $f(x)$, i.e., $D(f)=256(-a)^3=(-1)(16a)^2a$ is a non-square in $\FP{q}$, and that requires $a$ to be a non-square in $\FP{q}$.

Conversely, if $a$ is a non-square in $\FP{q}$, so is $D(f)$. By \cref{T51}, $f(x)$ is either irreducible over $\FP{q}$ or factored as a product of two linear factors and an irreducible quadratic factor over $\FP{q}$. The latter case cannot happen. Indeed, assume that
$$\begin{aligned}
x^4-a & =\left(x^2+rx+s\right)\left(x^2+tx+u\right) \\
      & =x^4+(r+t)x^3+(rt+s+u)x^2+(ru+st)x+su
\end{aligned}$$
where $r,s,t,u\in \FP{q}$. Matching coefficients of corresponding terms, we have that
$$\begin{aligned}
\begin{cases}
r+t    & =0  \\
rt+s+u & =0  \\
ru+st  & =0  \\
su     & =-a
\end{cases}
& \qquad \Rightarrow\qquad &
\begin{cases}
t       & =-r \\
s+u     & =r^2  \\
r(-s+u) & =0  \\
su      & =-a
\end{cases}
\end{aligned}$$
If $r=0$, then $u=-s$. So $a=-su=s^2$. If $r\neq 0$, then $u=s$. So $a=(-1)s^2$. In both cases, $a$ is a square in $\FP{q}$, which is impossible. Therefore $f(x)$ is irreducible. \mbox{}\hfill $\square$\\

Based on \cref{P52}, a class of stable quartic polynomials can be obtained.

\begin{thm}\label{T53}
Let $q$ be a prime power such that $q\equiv 1\Mod{4}$. If $a,b\in \FP{q}^\ast$, then
$$f(x)=b^3(x+a)^4-a$$
is stable over $\FP{q}$ iff $ab$ is a non-square in $\FP{q}$.
\end{thm}

\noindent \textit{Proof}\quad Assume that $q=p^m$ for some $m\in \NN$ and odd prime $p$. We only give the proof for the case where $p=8k+1$ since the other cases (where $p\equiv 3,5,7\Mod{8}$) can be proven in the same way.

If $f(x)$ is stable over $\FP{q}$, then it is irreducible over $\FP{q}$. Therefore, we know that $bf(x)=(b(x+a))^4-ab$ is irreducible over $\FP{q}$. The change of variables $y=b(x+a)$ then implies that $y^4-ab$ is irreducible over $\FP{q}$. From \cref{P52}, we deduce that $ab$ is a non-square in $\FP{q}$.

Conversely, assume that $ab$ is a non-square in $\FP{q}$. As before, following the idea introduced in \cref{sec2}, we construct a sequence $(\alpha_n)_{n\geq 0}$ such that $\alpha_0=0$ and that $f(\alpha_n)=\alpha_{n-1}\ (n\in \NN)$. For each $n\in \NI$, let
$$f_n(x)=f(x)-\alpha_n=b^3(x+a)^4-(\alpha_n+a).$$
By \cref{C28}, $f(x)$ is stable over $\FP{q}$ iff $f_n(x)$ is irreducible over $\FF{q}{4^n}$ for each $n\in \NI$. We see from the change of variables $y=b(x+a)$ that latter claim is true iff
$$g_n(y)=y^4-\beta_n$$
is irreducible over $\FF{q}{4^n}$ for each $n\in \NI$, where
$$\beta_n=b(\alpha_n+a)$$
By \cref{P52}, $f(x)$ is stable over $\FP{q}$ iff $\beta_n$ is a non-square in $\FF{q}{4^n}$ for each $n\in \NI$. It suffices to verify the latter claim for $n\in \NN$ because we know by assumption that $\beta_0=ab$ is a non-square in $\FP{q}$. Assume towards a contradiction that the claim is false. Then there exists a smallest $n_0\in \NN$ such that $\beta_{n_0}$ is a square in $\FF{q}{4^{n_0}}$.

If $n_0-1=0$, then $\beta_{n_0-1}=\beta_0$ is a non-square in $\FP{q}$. If $n_0-1\geq 1$, then by minimality of $n_0$, $\beta_{n_0-1}$ is again a non-square in $\FF{q}{4^{n_0-1}}$. By \cref{P52}, $g_{n_0-1}(x)$ is irreducible over $\FF{q}{4^{n_0-1}}$. If $x_0$ is a root of $g_{n_0-1}(x)$ in $\FF{q}{4^{n_0}}$, then $x_0^4=\beta_{n_0-1}$ and
$$\FF{q}{4^{n_0}}\cong \FF{q}{4^{n_0-1}}[x]\slash \left<g_{n_0-1}(x)\right>\cong \left\{r+sx_0+tx_0^2+ux_0^3:\ r,s,t,u\in \FF{q}{4^{n_0-1}}\right\}$$
By assumption, $\beta_{n_0}$ is a square in $\FF{q}{4^{n_0}}$. Hence, there exist $r,s,t,u\in \FF{q}{4^{n_0-1}}$ such that $\beta_{n_0}=\left(r+sx_0+tx_0^2+ux_0^3\right)^2$. Since $f(\alpha_n)=\alpha_{n-1}$ for each $n\in \NN$, we have that
\begin{equation}\label{Eq11}
\beta_n^4=\beta_{n-1}
\end{equation}
In particular, we know that
\begin{equation}\label{Eq12}
\left(r+sx_0+tx_0^2+ux_0^3\right)^8=\beta_{n_0}^4=\beta_{n_0-1}=x_0^4
\end{equation}
from which it follows that
\begin{equation}\label{Eq13}
\left(r+sx_0+tx_0^2+ux_0^3\right)^2=\epsilon x_0    
\end{equation}
where $\epsilon\in \FP{q}$ is a $4$-th root of unity. It is easy to see that at least three of $r,s,t,u$ must be non-zero. Indeed, if two of them, say $r$ and $s$ are $0$, then by \cref{Eq13},
$$t^2\beta_{n_0-1}+2tu\beta_{n_0-1}x_0+u^2\beta_{n_0-1}x_0^2=\epsilon x_0$$
So $t^2\beta_{n_0-1}=u^2\beta_{n_0-1}=0$. However, by \cref{Eq11}, $\beta_n\neq 0$ for each $n\in \NN$ since $\beta_0\neq 0$. Thus, $t=u=0$, meaning that $\beta_{n_0}=0$, which is impossible. The other cases can be proven in the same way.

Meanwhile, recall that $p=8k+1$. So \cref{Eq12} implies that
$$\left(r+sx_0+tx_0^2+ux_0^3\right)^{8k+1}=\beta_{n_0-1}^k\left(r+sx_0+tx_0^2+ux_0^3\right)$$
which, after simplification, becomes
$$r^p+s^p\beta_{n_0-1}^{2k}x_0+t^p\beta_{n_0-1}^{4k}x_0^2+u^p\beta_{n_0-1}^{6k}x_0^3=\beta_{n_0-1}^k\left(r+sx_0+tx_0^2+ux_0^3\right)$$
Matching coefficients of corresponding terms, we then obtain that
\begin{equation}\label{Eq14}
\begin{cases}
r^p & = r\beta_{n_0-1}^k \\
s^p\beta_{n_0-1}^{2k} & = s\beta_{n_0-1}^k
\end{cases}
\end{equation}
Previous discussion indicates that one of $r$ and $s$ must be non-zero. Without loss of generality, we may consider only the case where $r\neq 0$ since the same proof remains valid in the other case. Since $r\neq 0$, $\beta_{n_0-1}^k=r^{p-1}=r^{8k}$. Thus, $\beta_{n_0-1}=vr^8$ for some $k$-th root of unity $v\in \overline{\FP{q}}$. Since $2k\mid (8k+1)^m-1=q-1$, $\FP{q}$ contains all the $2k$-th roots of unity (and hence all the $k$-th roots of unity), say $\pm v_j\ (1\leq j\leq k)$, where $v_{j_1}\neq \pm v_{j_2}$ whenever $j_1\neq j_2$. Thus, all $k$-th roots of unity are squares in $\FP{q}$ because they are simply $v_j^2\ (1\leq j\leq k)$. In particular, $\beta_{n_0-1}$ is a square in $\FF{q}{4^{n_0-1}}$, which, as has been shown earlier, is impossible. Hence, no such $n_0$ exists, proving that $f(x)$ is stable over $\FP{q}$.\mbox{}\hfill $\square$\\

\par Next, let $q=3^m$. We study the stability of polynomials over $\FP{q}$ of the form
$$f(x)=b^3(x+a)^4+c(x+a)-a$$
where $a,b,c\in \FP{q}$, and $b,c\neq 0$. As before, we construct a sequence $(\alpha_n)_{n\geq 0}$ such that $\alpha_0=0$ and that $f(\alpha_n)=\alpha_{n-1}\ (n\in \NN)$. For each $n\in \NI$, let
$$f_n(x)=f(x)-\alpha_n=b^3(x+a)^4+c(x+a)-(\alpha_n+a)$$
By \cref{C28}, $f(x)$ is stable over $\FP{q}$ iff $f_n(x)$ is irreducible over $\FF{q}{4^n}$ for each $n\in \NI$. The change of variables $y=b(x+a)$ implies that latter claim is true iff
$$g_n(y)=y^4+cy-\beta_n$$
is irreducible over $\FF{q}{4^n}$ for each $n\in \NI$, where
$$\beta_n=b(\alpha_n+a)$$ 
To proceed, we now prove a theorem describing the factorization patterns of the polynomials $g_n(x)$ over $\FP{3^m}$.

\begin{thm}\label{T54}
Let $q=3^m$ and $c\in \FP{q}^\ast$.  If $g(x)=x^3-dx-c^2$ and $h(x)=x^4+cx+d$, then the following hold.
\begin{enumerate}
    \item $h(x)$ splits completely into linear factors over $\FP{q}$ iff $g\left(r_j^2\right)=0\ (1\leq j\leq 3)$ for some pairwise distinct $r_1,r_2,r_3\in \FP{q}$ such that $-r_j^2\pm cr_j^{-1}\ (1\leq j\leq 3)$ are squares in $\FP{q}$; \label{I1} 
    \item $h(x)$ is factored as a product of two irreducible quadratic factors over $\FP{q}$ iff $g\left(u_j\right)=0\ (1\leq j\leq 3)$ for some pairwise distinct $u_1,u_2,u_3\in \FP{q}$ such that $u_1,u_2$ are non-squares in $\FP{q}$ and that $u_3=r^2$ for some $r\in \FP{q}$, where $-r^2\pm cr^{-1}$ are non-squares in $\FP{q}$; \label{I2}
    \item $h(x)$ has a unique root $\FP{q}$ iff $g(x)$ is irreducible over $\FP{q}$; \label{I3}
    \item $h(x)$ has exactly two roots in $\FP{q}$ iff $g(x)$ has a unique root $r^2\in \FP{q}$, where $r\in \FP{q}$ and $r^4-c^2r^{-2}$ is a non-square in $\FP{q}$; \label{I4}
    \item $h(x)$ is irreducible over $\FP{q}$ iff $g(x)$ has a unique root which is a non-square in $\FP{q}$. \label{I5}
\end{enumerate}
\end{thm}

\noindent\textit{Proof}\quad First, we observe that if $h(x)$ is factored into a product of two (possibly reducible) quadratic factors over $\FP{q}$, say
$$\begin{aligned}
x^4+cx+d & =\left(x^2+rx+s\right)\left(x^2+tx+u\right) \\
        & =x^4+(r+t)x^3+(rt+s+u)x^2+(ru+st)x+su
\end{aligned}$$
Then we have that
$$\begin{aligned}
\begin{cases}
r+t    & =0 \\
rt+s+u & =0 \\
ru+st  & =c \\
su     & =d
\end{cases}
& \qquad \Rightarrow\qquad &
\begin{cases}
t       & =-r \\
s+u     & =r^2 \\
r(-s+u) & =c \\
su      & =d
\end{cases}
\end{aligned}$$
By the third equation, $r\neq 0$. Solving the second and third equation for $s,u$ yields that $s=cr^{-1}-r^2,u=-cr^{-1}-r^2$. These, together with the fourth equation, imply that $r^4-c^2r^{-2}=su=d$. Thus, $g\left(r^2\right)=r^6-dr^2-c^2=0$. So $h(x)$ can be factored as a product of two quadratic factors over $\FP{q}$ iff there is an $r\in \FP{q}$ such that
\begin{equation}\label{Eq15}
f(x)=\left(x^2+rx+cr^{-1}-r^2\right)\left(x^2-rx-cr^{-1}-r^2\right)
\end{equation}
and the latter claim is true iff $g(x)$ has a root which is a square in $\FP{q}$. We also note that the discriminant of the first quadratic factor in \cref{Eq15} is $-r^2-cr^{-1}$ and that of the second quadratic factor is $-r^2+cr^{-1}$.\\

\par First, we prove item (\ref{I1}). Clearly, $h(x)$ is separable over $\FP{q}$. Assume that $h(x)$ has four pairwise distinct roots. Then it is easy to check that they give rise to exactly three pairwise distinct factorizations in the form of \cref{Eq15}, and in each factorization, both quadratic factors are reducible over $\FP{q}$. So Item (\ref{I1}) follows from the remarks on \cref{Eq15}.

Next, we prove Item (\ref{I2}). If $h(x)$ is factored as a product of two irreducible quadratic factors, then there is an $r\in \FP{q}$ satisfying \cref{Eq15} such that $r^2\in \FP{q}$ is a root of $g(x)$ and that $-r^2\pm cr^{-1}$ are non-squares in $\FP{q}$. Also, since $h(x)$ has exactly two distinct irreducible factors over $\FP{q}$, \cref{T51} implies that the discriminant of $h(x)$, i.e., $D(h)=256d^3$ is a square in $\FP{q}$. Thus, $d$ is a square in $\FP{q}$. By \cite[Theorem 2]{KW75a}, $g(x)$ either has three roots in $\FP{q}$ or is irreducible over $\FP{q}$. The latter case cannot happen since we have shown that $g(x)$ has a root in $\FP{q}$. Hence, $g(x)$ has two other roots in $\FP{q}$ which must be non-squares in $\FP{q}$ (because if one of them is a square, then \cref{Eq15} gives rise to another factorization of $h(x)$ into irreducible factors over $\FP{q}$). The other direction is straightforward.

Then we prove Item (\ref{I3}). If $h(x)$ is factored as a product of a linear factor and an irreducible cubic factor, then like in the previous proof, \cref{T51} and \cite[Theorem 2]{KW75a} imply that $g(x)$ either has three roots $u_1,u_2,u_3$ in $\FP{q}$ or is irreducible over $\FP{q}$. In the former case, we know from Vieta's formula that $u_1u_2u_3=c^2$. If all of $u_1,u_2,u_3$ are non-squares in $\FP{q}$, then so is $c^2$, which is impossible. So one of them is a square in $\FP{q}$, giving rise to a factorization of $h(x)$ into two quadratic factors via \cref{Eq15}. Since that is impossible, $g$ must be irreducible over $\FP{q}$. The other direction is straightforward.

Now we prove Item (\ref{I4}). If $h(x)$ is factored as a product of two linear factors and an irreducible quadratic factor, then there is an $r\in \FP{q}$ satisfying \cref{Eq15} such that $r^2$ is a root of $g(x)$ and that exactly one of $-r^2\pm cr^{-1}$ is a square in $\FP{q}$. Meanwhile, \cref{T51} implies that the discriminant of $h(x)$ is a non-square in $\FP{q}$. Thus, $d$ is a non-square in $\FP{q}$ and by \cite[Theorem 2]{KW75a}, $g(x)$ has no other root in $\FP{q}$. The other direction is straightforward.

Finally, we prove Item (\ref{I5}). If $h(x)$ is irreducible over $\FP{q}$, then like in the previous proof, \cref{T51} and \cite[Theorem 2]{KW75a} imply that $g(x)$ has a unique root $u$ in $\FP{q}$. Since no $r\in \FP{q}$ satisfies \cref{Eq15}, $u$ is a non-square in $\FP{q}$.\mbox{}\hfill $\square$\\

\par Let $a,b,c$ and $(\beta_n)_{n\geq 0}$ be as previously defined. Then \cref{T54} leads to the following result.

\begin{Cor}\label{C55}
The polynomial $f(x)=b^3(x+a)^4+c(x+a)-a$ is stable over $\FF{3}{m}$ iff for each $n\in \NI$, $x^3-\beta_nx-c^2$ has a unique root which is a non-square in $\FF{3}{4^nm}$.
\end{Cor}

\section{Conclusions}\label{sec6}
In this paper, several classes of stable shift-resistant polynomials of degrees $2,3$ and $4$ are constructed by repeated use of Capelli's lemma.  As a result, we confirm 
the conjecture of Ahmadi and Monsef-Shokri \cite{AMS20a}  that the polynomial $f(x) = x^3 + x^2 + 1$ is stable over $\FP{2}$. 
When $q$ is a prime power and $q\equiv 1\Mod{4}$, an improvement to the lower bound of the number of stable quadratic polynomials over $\FP{q}$ is obtained. No lower bound of the number of stable polynomials of an arbitrary degree $d$ over finite fields seems to have been established (although upper bounds have been studied in \cite[Theorem 1]{GPN10a} and \cite[Theorem 5.5]{GPNOS14a} using character sums). For the future study,   we are hopeful that based on the stability test in \cref{sec2}, stable quadratic polynomials over $\FP{q}$ can be constructed when $q\equiv 3\Mod{4}$ and examples of stable polynomials of degree $\geq 5$ can be found. Also, we are working towards establishing a better estimate on the number of stable polynomials over finite fields.

\setcounter{secnumdepth}{0}

\end{document}